\documentclass[12pt]{article}
\usepackage[latin1]{inputenc}
\usepackage[activeacute,spanish]{babel}
\usepackage{amssymb,amsmath}
\usepackage{epsfig}
\renewcommand{\aa}{\ensuremath{\mathbf{a}}}

\newcommand{\abs}[1]{\ensuremath{\left|#1\right|}}

\newcommand{\Div}{\nabla\cdot}
\newcommand{\ee}{\ensuremath{\mathbf{e}}}
\newcommand{\eps}{\ensuremath{\varepsilon}}
\newcommand\escala{\tilde{\varphi}}

\newcommand{\Grad}{\ensuremath{\mathrm{\nabla}}}

\newcommand{\M}{\ensuremath{\mathcal{M}}}
\newcommand{\N}{\ensuremath{\mathbb{N}}}
\newcommand{\n}{\ensuremath{\mathbf{n}}}
\newcommand{\Om}{\ensuremath{\Omega}}
\newcommand{\pt}{\ensuremath{\partial_t}}

\renewcommand{\r}{\ensuremath{\mathbf{r}}}
\newcommand{\RR}{\mathbb{R}}

\newcommand\wavelet{\tilde{\psi}}
\newcommand{\x}{\ensuremath{\mathbf{x}}}

\renewcommand{\u}{\mathbf{u}}

\newcommand{\Iap}{I_{\mathrm{app}}}
\newcommand{\Ion}{I_{\mathrm{ion}}}

\newcommand{\IapK}{I_{\mathrm{app,K}}}

\newcommand{\IonK}{I_{\mathrm{ion},K}}

\newcommand{\dx}{\ensuremath{\, dx}}

\newcommand{\EE}{\mathcal{E}}
\newcommand{\TT}{\mathcal{T}}

\newcommand{\bM}{\ensuremath{\mathbf{M}}}
\begin{document}
\thispagestyle{empty}

\begin{minipage}{0.4\textwidth}
\includegraphics[width=0.3\textwidth]{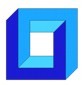}

{\bf Sociedad Chilena de\\
Mecánica Computacional}
\end{minipage}
\begin{minipage}{0.5\textwidth}
\begin{center}
{\bf Cuadernos de Mecánica Computacional\\

\vspace{.2cm}
Vol. 6 nº1, 2008}
\end{center}
\end{minipage}

\vspace{.6cm}
\begin{center}
{\bf UN MÉTODO ADAPTATIVO PARA EL MODELO BIDOMINIO EN ELECTROCARDIOLOGÍA \\

\vspace{.5cm}
Mostafa Bendahmane, Raimund Bürger y Ricardo Ruiz Baier$^*$}\\

\vspace{.5cm}
{\small $^*$ Departamento de Ingeniería Matemática - Universidad de Concepción\\
Casilla 160-C - Concepción - CHILE\\
e-mail : {\tt mostafab,rburger,rruiz@ing-mat.udec.cl}}\\
\end{center}

\begin{abstract}
En este trabajo se presenta un método de volúmenes finitos enriquecido con un
esquema de multiresolución completamente adaptativo para obtener adaptatividad espacial,
y un esquema Runge-Kutta-Fehlberg con paso temporal de variación local para obtener
adaptatividad temporal, para resolver numéricamente las conocidas ecuaciones "bidominio"
que modelan la actividad eléctrica del tejido en el miocardio. Se consideran dos modelos
simples para las corrientes de membrana y corrientes iónicas. En primer lugar definimos
 una solución aproximada y nos referimos a su convergencia  a la correspondiente
solución débil del problema continuo, obteniendo de este modo una demostración alternativa
de que el problema continuo es bien puesto. Luego de introducir la técnica de multiresolución,
se deriva un umbral óptimo para descartar la información no significativa, y tanto la
eficiencia como la precisión de nuestro método es vista en términos de la aceleración de
tiempo de máquina, compresión de memoria computacional y errores en diferentes normas.
\end{abstract}
\section{Introducción}
Las mediciones directas representan una dificultad obvia en ciencias. Por lo tanto
simulaciones numéricas son de gran interés, específicamente en modelos cardíacos.
Entre tales modelos, el \emph{modelo bidominio} es conocido como uno de los más
precisos y completos para el estudio teórico y numérico de la actividad eléctrica
en el tejido cardíaco. Desde el punto de vista computacional, el modelo bidominio
representa un verdadero
desafío, dado que el tejido cardíaco tiene tamaños del orden de centímetros y por
ejemplo, los frentes de exitación de las ondas son del orden de los $10^{-3}$
centímetros. Esta característica local, no sólo espacial, sino también temporal,
junto con la aparición de frentes perfilados en el campo de los potenciales eléctricos,
hace que las simulaciones en mallas uniformes sean prácticamente imposibles de llevar
a cabo. Aquí es donde los métodos adaptativos juegan un rol vital en las simulaciones
cardíacas. En este artículo desarrollamos un esquema de multiresolución completamente adaptativo
provisto de adaptatividad temporal a traves de una estrategia de paso temporal local
y deducimos un umbral óptimo para descartar información no significativa. Basados en
experiencia previa  sobre sistemas de reacción-difusión y ecuaciones
parabólicas degeneradas \cite{bbrs,brss,brss2}, sugerimos que la multiresolución puede
ser una herramienta efectiva para resolver las ecuaciones del modelo bidominio.
En el contexto de métodos de multiresolución completamente adaptativos, mencionamos que
existen trabajos desarrollados por varios grupos de investigación (ver
\cite{bbrs,brss,Cohen,Muller,MS}), con aplicaciones enfocadas a otras áreas.
\subsection{El modelo bidominio}
Supongamos $\Om\subset \RR^2$ abierto y acotado con frontera suave $\partial \Om$.
En este modelo, el músculo cardíaco es interpretado como la unión de dos medios
continuos interpenetrados y superimpuestos: el medio intracelular y el extracelular.
Estos ocupan la misma área y se encuentran separados por la membrana celular cardíaca.
$u_\mathrm{i}=u_\mathrm{i}(t,x)$ y $u_\mathrm{e}=u_\mathrm{e}(t,x)$ representan
los potenciales eléctricos \textit{intracelular} y \textit{extracelular} en
$(x,t)\in \Om_T:=\Om\times(0,T)$, y la diferencia entre estos potenciales
$v=v(t,x)=u_\mathrm{i}-u_\mathrm{e}$ es conocido como el potencial
\textit{transmembrana}. La conductividad del tejido está representada por
tensores escalados  $\bM_\mathrm{i}(x)$ y $\bM_\mathrm{e}(x)$ dados por
$$\bM_\mathrm{j}(x)=\sigma_\mathrm{j}^t\mathbf{I} +
(\sigma_\mathrm{j}^l-\sigma_\mathrm{j}^t)\aa_l(x)\aa_l^T(x),$$
donde $\sigma_\mathrm{j}^l=\sigma_\mathrm{j}^l(x)\in C^1(\RR^2)$ y
$\sigma_\mathrm{j}^t=\sigma_\mathrm{j}^t(x)\in C^1(\RR^2)$, para
$\mathrm{j}=\mathrm{i,e}$  son las conductividades intra- y extracelulares a lo
largo y a través respectivamente de la dirección de la fibra muscular correspondiente
(paralela a $\aa_l(x)$).

Comúnmente se utilizan los radios de anisotropía $\frac{\sigma_\mathrm{i}^l}{\sigma_\mathrm{i}^t}$
y $\frac{\sigma_\mathrm{e}^l}{\sigma_\mathrm{e}^t}$. En general las conductividades
en la dirección longitudinal $l$ son de mayor magnitud que aquellas a través de la
fibra (dirección $t$); y tal caso se denomina \emph{anisotropía fuerte} en conductividad
eléctrica.

El siguiente sistema fuertemente acoplado de reacción-difusión forma el modelo bidominio
anisotrópico (ver \cite{ying}):
\begin{equation}\label{S4}
\begin{split}
\beta c_m\pt v+\Div(\bM_\mathrm{e}(x)\Grad u_\mathrm{e})+\beta \Ion(v,w)&=\Iap \textrm{ en }\Om_T,\\
-\Div ((\bM_\mathrm{i}(x)+\bM_\mathrm{e}(x))\Grad u_\mathrm{e})-
\Div(\bM_\mathrm{i}(x)\Grad v)&=0 \textrm{ en }\Om_T,\\
\pt w -H(v,w)&=0 \textrm{ en }\Om_T.
\end{split}
\end{equation}
Aquí, $c_m>0$ representa la capacitancia de superficie de la membrana, $\beta$ es la
razón área-volumen, y $w(t,x)$ es la variable de recuperación, que toma en cuenta
las variables de concentración del modelo. Las corrientes de los estímulos aplicados
a los medios intra- y extracelulares están representadas por la función $\Iap=\Iap(t,x)$
que satisface $\int_\Om \Iap(t,x)\dx=0$ para casi todo $t\in (0,T)$. La función $H$ en
la ecuación diferencial ordinaria de (\ref{S4}) y la función $\Ion$
corresponden a uno de los modelos más simples para las corrientes de la membrana y
corriente iónica (entre una amplia variedad de tales modelos): el modelo de membrana
de Mitchell--Schaeffer \cite{MS:Ion}
\begin{equation*}
H(v,w)=\frac{w_\infty(v/v_p)-w}{R_mc_m\eta_\infty(v/v_p)},\quad
\Ion(v,w)=\frac{v_p}{R_m}\Biggl(\frac{v}{v_p\eta_2}-
\frac{v^2(1-v/v_p)w}{v_p^2\eta_1}\Biggr),
\end{equation*}
donde
$$\eta_\infty(s)=\left\{\begin{array}{cc}
\eta_3& \textrm{ para } s<\eta_5,\\
\eta_4& \textrm{ en otro caso},\end{array}\right.,\quad
w_\infty(s)=\left\{\begin{array}{cc}
1& \textrm{ para } s<\eta_5,\\
0&  \textrm{ en otro caso}.\end{array}\right. $$
$R_m$ representa la resistividad superficial de la membrana y
$v_p,\eta_1,\eta_2,\eta_3,\eta_4,\eta_5$ son constantes dadas. El sistema
(\ref{S4}) es provisto con condiciones de borde de
no-flujo, representando un tejido cardíaco aislado
\begin{equation}\label{S2}
(\bM_\mathrm{j}(x)\Grad u_\mathrm{j})\cdot\n =0\textrm{ sobre }
\Sigma_T:=\partial\Om\times(0,T),\quad \mathrm{j}=\mathrm{i,e},
\end{equation}
y condiciones iniciales apropiadas en $\Om$ para el potencial
transmembrana y variable de recuperación $v(0,x)=v_0(x)$, $w(0,x)=w_0(x)$.
Para asegurar la dependencia continua de los datos en la componente $v$
de la solución, requerimos que el dato inicial $v_0$ sea compatible con
(\ref{S2}). Por lo tanto la condición de compatibilidad
\begin{equation}\label{compat_ue}
\int_\Om u_\mathrm{e}(x,t)\dx=0 \textrm{ para c.t. }t\in(0,T),
\end{equation}
debe ser satisfecha.

La teoría estándar para ecuaciones parabólicas--elípticas no puede ser aplicada
de forma natural en el análisis de las ecuaciones del modelo bidominio, debido a
la diferencia existente entre los grados de anisotropía entre los medios intra- y
extra celulares. Debido a esta característica, el sistema  (\ref{S4}) es de
naturaleza parabólica degenerada. En \cite{Bend-Karl:cardiac} los autores prueban
existencia y unicidad de solución para las ecuaciones de bidominio, utilizando
el método de Faedo--Galerkin y teoría de compacidad.

\section{Un método base de volúmenes finitos}
Para resolver numéricamente (\ref{S4}) introducimos un método
estándar de volúmenes finitos. Una malla admisible para $\Om$ será  formada
por una familia $\TT$ de volúmenes de control (polígonos abiertos y convexos)
de diámetro máximo $h$. Para todo $K \in \TT$, $x_K$ denota el centro de $K$,
$N(K)$ el conjunto de vecinos de $K$, $\EE_{\textrm{int}}(K)$ el conjunto de
bordes de $K$ en el interior de $\TT$ y $\EE_{\textrm{ext}}(K)$
el conjunto de bordes de $K$ sobre la frontera $\partial \Om$. Para todo $L \in N(K)$
$d(K,L)$ denota la distancia entre $x_K$ y $x_{L}$, $\sigma_{K,L}$ es
la interfaz entre $K$ y $L$, y $\eta_{K,L}$ ($\eta_{K,\sigma}$
respectivamente) es el vector unitario normal a $\sigma_{K,L}$
orientado desde $K$ hacia $L$. Para todo $K \in \TT$, $|K|$
es la medida de $K$. La admisibilidad de $\TT$ implica que
$\overline{\Om}=\cup_{K\in \TT} \overline{K}$, $K\cap L=\emptyset$
si $K,L\in \TT$ y $K \ne L$, y además existe una sucesión finita
$(x_{K})_{K\in \TT}$, tal que $\overline{x_{K}x_{L}}$ es ortogonal a
$\sigma_{K,L}$. Ahora, considerar $K \in \TT$ y $L \in N(K)$ con vértices comunes
$(a_{\ell,K,L})_{1\le \ell\le I}$ con $I \in \N \backslash \{0\}$, y
denotemos por $T_{K,L}$ al polígono abierto y convexo de vértices
$(x_K,x_L)$ y $(a_{\ell,K,L})_{1\le \ell\le I}$.
Sea $\mathcal{D}$ una discretización admisible de $Q_T$, que consiste en una
malla admisible para $\Om$, un paso temporal $\Delta t>0$, y $N>0$ elegido
como el menor entero tal que $N\Delta t\ge T$. Con esto, escribimos
$t^n=n \Delta t$ para $n\in  \{0,\ldots,N\}$. Sobre cada elemento $K \in \TT$,
se definen tensores de conductividad (definidos positivos) mediante
$$M_{\mathrm{j},K}=\frac{1}{|K|}\int_\Om \bM_{\mathrm{j}}(x)\dx,\quad
\mathrm{j}=\mathrm{i,e}.$$
Sea $F_{j,K,L}$ una aproximación de $\int_{\sigma_{K,L}}\bM_\mathrm{j}(x)\Grad u_\mathrm{j} \cdot \eta_{K,L}d\gamma$
para $\mathrm{j}=\mathrm{i,e}$, y para $K \in \Om_R$, $L \in N(K)$, sea
$$M_{\mathrm{j},K,L}=\abs{\frac{1}{|K|} \int_K \bM_\mathrm{j}(x)\dx \,\eta_{K,L}}
:=\abs{M_{\mathrm{j},K}\,\eta_{K,L}},\quad \mathrm{j}=\mathrm{i,e}.$$
Los flujos difusivos $M_\mathrm{j}(x)\Grad u_j \cdot \eta_{K,L}$ sobre $\sigma_{K,L}$
son aproximados por
$$ \int_{\sigma_{K,L}} (\bM_\mathrm{j}(x)\Grad u_j) \cdot \eta_{K,L} d\gamma
\approx |\sigma_{K,L}|M_{\mathrm{j},K,L}
\frac{u_{\mathrm{j},\sigma}-u_{\mathrm{j},K}}{d(K,\sigma_{K,L})},$$
donde $y_\sigma$ es el centro de $\sigma_{K,L}$ y $u_{\mathrm{j},\sigma}$ es una
aproximación de $u_\mathrm{j}(y_\sigma)$, $\mathrm{j}=\mathrm{i,e}$. La conservatividad
del método nos permite determinar las incógnitas adicionales $u_{\mathrm{j},\sigma}$,
y además calcular los flujos numéricos sobre los bordes:
\begin{equation*}
\textrm{$F_{j,K,L}=d_{j,K,L}^*
\frac{|\sigma_{K,L}|}{d(K,L)}(u_{\mathrm{j},L}-u_{\mathrm{j},K})$ si $L\in N(K)$},
\end{equation*}
donde
$$d_{\mathrm{j},K,L}^*=\frac{M_{\mathrm{j},K,L}M_{\mathrm{j},L,K}}
{d(K,\sigma_{K,L})M_{\mathrm{j},K,L}+d(L,\sigma_{K,L})M_{\mathrm{j},L,K}} d(K,L).
$$
Finalmente, aproximaremos el sistema (\ref{S4}) mediante la siguiente formulación de
volúmenes finitos: Determinar $(u_{\mathrm{j},K}^{n})_{K \in \TT}$ para $\mathrm{j}=\mathrm{i,e}$
y $n\in \{0,\ldots,N\}$,  $(v_{K}^{n})_{K \in \TT}
=(u_{\mathrm{i},K}^{n}-u_{\mathrm{e},K}^{n})_{K \in \TT}$ para $n\in \{0,\ldots,N\}$, y
$(w_{K}^{n})_{K \in \TT}$ para $n\in \{0,\ldots,N\}$, tales que para todo
$K \in \TT$ y $n \in  \{0,\ldots,N-1\}$
\begin{equation}\label{prob:init}
v_K^0=\frac{1}{|K|} \int_{K} v_0(x) \dx,\quad
w_K^0=\frac{1}{|K|} \int_{K} w_0(x) \dx,
\end{equation}
\begin{equation*}
\begin{split}
\beta c_{\mathrm{m}} |K|\frac{v^{n+1}_K-v^{n}_K}{\Delta t}+\sum_{L \in N(K) }
     d_{\mathrm{e},K,L}^*\frac{|\sigma_{K,L}|}{d(K,L)}
(u^{n}_{\mathrm{e},L}-u^{n}_{\mathrm{e},K})
     +\beta |K|\IonK^n
& =|K|\IapK^{n},
\\
\sum_{L \in N(K) } \frac{|\sigma_{K,L}|}{d(K,L)}    \left\{
\bigl(d_{\mathrm{i},K,L}^*+d_{\mathrm{e},K,L}^*\bigr)
\bigl(u^{n+1}_{\mathrm{e},L}-u^{n+1}_{\mathrm{e},K}\bigr)
+ 
 d_{i,K,L}^*\bigl(v^{n+1}_L-v^{n+1}_K\bigr)\right\} &=|K|\IapK^{n},
\\
     |K|\frac{w^{n+1}_K-w^{n}_K}{\Delta t}-|K|H^{n}_K&=0.
\end{split}
\end{equation*}

La condición de borde (\ref{S2}) es tomada en cuenta imponiendo condiciones de
no-flujo sobre los bordes externos:
\begin{equation}\label{moyenne-discr}
d_{\mathrm{j},K,\sigma}^*\frac{|\sigma_{K,L}|}{d(K,L)}(u^{n}_{\mathrm{j},L}-u^{n}_{\mathrm{j},K})=0
\textrm{ for }\sigma \in \EE_{\textrm{ext}}(K),\quad \mathrm{j}=\mathrm{i,e},
\end{equation}
y (\ref{compat_ue}) es discretizada mediante $\sum_{K\in \TT}|K|u^n_{\mathrm{e},K}=0$,
para todo $n\in \{1,\ldots,N-1\}$.

La existencia, unicidad de solución aproximada y convergencia del esquema numérico hacia la solución
débil correspondiente,  son analizados en el trabajo \cite{bbrconvergence}. Aún más, como en
\cite{bbrs}, es posible deducir que en el caso de mallas uniformes, el esquema
(\ref{prob:init})-(\ref{moyenne-discr}) es estable bajo la condición CFL
\begin{equation}\label{cfl}
\Delta t\leqslant\frac{h}{2\displaystyle{\max_{K\in\TT}\Bigl(|\IonK|+2|\IapK|\Bigr)}+
4h^{-1}\max_{K\in\TT}\Bigl(|M_{i,K}|+|M_{e,K}|\Bigr)}.
\end{equation}
\section{Multiresolución y wavelets}
Considerar como dominio computacional, un simple rectángulo que luego de un cambio de
variables corresponde a $\Om=[0,1]^2$. En primer lugar introducimos una jerarquía de
mallas anidadas \mbox{$\Lambda_0\subset\cdots\subset\Lambda_L$}, usando una partición
diádica uniforme de $\Om$. Cada malla $\Lambda_l:=\{V_{(i,j),l}\}_{(i,j)}$, con
$(i,j)$ a ser definido, está formada por volúmenes de control en cada nivel de resolución
$V_{(i,j),l}:=2^{-l}[i,i+1]\times[j,j+1]$, $i,j\in I_l=\{0,\ldots,2^l-1\}$,
$l=0,\ldots,L$. $l=0$ corresponde al nivel más grueso y $l=L$ al más fino.
Introducimos también los conjuntos de refinamiento  $\M_{(i,j),l}=\{2(i,j)+\ee\}$,
$\ee\in E:=\{0,1\}^2$ con $\#\M_{(i,j),l}=4.$ Para cada nivel $l=0,\ldots,L$, se define
la función de escala
\begin{equation*}
\escala_{(i,j),l}(\x)=\frac{1}{|V_{(i,j),l}|}\chi_{V_{(i,j),l}}(\x)=
2^{2l}\chi_{[0,1]^2}(2^lx_1-i,2^lx_2-j),
\end{equation*}
y por lo tanto el promedio de $u(\cdot,t)\in L^1(\Om)$ sobre el volumen de control
$V_{(i,j),l}$ puede expresarse en términos del producto interior
$\bar{u}_{(i,j),l}:=\langle u,\escala_{(i,j),l}\rangle_{L^1(\Om)}.$
Con esto en mente, es posible definir una relación de dos escalas para
las funciones de escala y medias en celda respectivamente
$$\escala_{(i,j),l}=\sum_{\r\in\M_{(i,j),l}}\frac{|V_{\r,l+1}|}{|V_{(i,j),l}|}
\escala_{\r,l+1},\qquad
\bar{u}_{(i,j),l}=\sum_{\r\in\M_{(i,j),l}}\frac{|V_{\r,l+1}|}{|V_{(i,j),l}|}\bar{u}_{\r,l+1}.$$
Tal relación define un operador de \emph{proyección} que transforma elementos desde niveles
finos a niveles gruesos. Para $\x\in V_{2(i,j)+\aa,l+1}$, $\aa\in E$,
definimos la función \emph{wavelet} en función de la función de escala sobre un nivel
más fino, como
\begin{align*}
\wavelet_{(i,j),\ee,l}&=\sum_{\aa\in E}2^{-2}(-1)^{\aa\cdot\ee}
\escala_{2(i,j)+\aa,l+1}
=\sum_{\r\in\M_{(i,j),l}}
\frac{|V_{\r=2(i,j)+\aa,l+1}|}{|V_{(i,j),l}|}(-1)^{\aa\cdot\ee}
\escala_{\r,l+1}.
\end{align*}
Por otro lado, para todo $\ee\in E^*:=E\setminus\{(0,0)\}$, es posible obtener una relación
de dos escalas \emph{inversa} (ver \cite{Muller})
$$\escala_{2(i,j)+\aa,l+1} =\sum_{\ee\in E}(-1)^{\aa\cdot\ee}
\wavelet_{(i,j),\ee,l} ,\quad \aa\in E.$$
Ahora, para $\ee\in E^*$, introducimos los \emph{detalles}, que juegan un papel
crucial en la detección de zonas donde la solución posee altos gradientes
\begin{equation*}
d_{(i,j),\ee,l}:=\langle u,\wavelet_{(i,j),\ee,l}\rangle.
\end{equation*}
Dada la relación de dos escalas inversa, es posible escribir:
\begin{equation}\label{tm_a}
\hat{u}_{(i,j),l+1}=\sum_{\r\in \bar{S}^l_{(i,j)}}g_{(i,j),\r}^l\bar{u}_{\r,l},
\end{equation}
donde $\bar{S}^l_{(i,j)}:=\bigr\{V_{([i/2]+r_1,[j/2]+
r_2),l}\bigl\}_{r_1,r_2\in\{-s,\ldots,0,\ldots,s\}}$ denota el  esténcil de interpolación,
$g_{(i,j),\r}^l$ son coeficientes de la interpolación, y el gorro sobre $u$ en el lado
izquierdo de (\ref{tm_a}) indica que se trata de un valor predicho. La relación (\ref{tm_a})
define un operador de \emph{predicción}, que transforma elementos de niveles  gruesos en niveles
finos. En contraste con el operador de proyección, el operador de predicción no es único. Sin
embargo, se imponen dos restricciones básicas para su definición: Debe ser consistente
con la proyección, en el sentido de que es el \emph{inverso a derecha} del operador de
proyección; y debe ser local, en el sentido de que el valor predicho dependerá sólo de
$\bar{S}^l_{(i,j)}$. En nuestro caso particular, la predicción es polinomial:
 \begin{equation*}
 \hat{u}_{2i+e_1,2j+e_2,l+1}=\bar{u}_{i,j,l}-(-1)^{e_1}Q_x-(-1)^{e_2}Q_y
 +(-1)^{e_1e_2}Q_{xy},
 \end{equation*}
 donde $e_1,e_2\in \{ 0,1\}$ y
 \begin{eqnarray*}
 Q_x&:=&\sum_{n=1}^s\tilde{\gamma}_n(\bar{u}_{i+n,j,l}-\bar{u}_{i-n,j,l}),
 Q_y:=\sum_{p=1}^s\tilde{\gamma}_p(\bar{u}_{i,j+p,l}-\bar{u}_{i,j-p,l}), \\
 Q_{xy}&:=&\sum_{n=1}^s\tilde{\gamma}_n\sum_{p=1}^s\tilde{\gamma}_p(\bar{u}_{i+n,j+p,l}-      \bar{u}_{i+n,j-p,l}-\bar{u}_{i-n,j+p,l}+\bar{u}_{i-n,j-p,l}).
 \end{eqnarray*}
Los coeficientes correspondientes son $\tilde{\gamma}_1=\smash{-\frac{22}{128}}$ y
$\tilde{\gamma}_2=\smash{\frac{3}{128}}$.
Cuanto más regular es la función $u$ sobre $V_{(i,j),l}$, el coeficiente de detalle correspondiente
es más pequeño en módulo. En vista de esta propiedad de cancelación, es natural pensar en
alguna estrategia para eliminar información no significativa (\emph{estrategia de corte}). La idea
básica es eliminar todos los
elementos de la malla que correspondan a detalles que se encuentran bajo una tolerancia (dependiente
del nivel de resolución) dada por
\begin{equation} \label{eq4.4}
 \eps_l=2^{2(l-L)}\eps_R,
\end{equation}
donde $\epsilon_R$ es una tolerancia de referencia a ser determinada en la sección \ref{sec:err-analysis}.
\subsection{Estructura de datos en árbol}
Organizaremos las medias en celda y los detalles correspondientes
utilizando una estructura de árbol graduado dinámico. Este tipo de
almacenamiento garantiza la estabilidad de las operaciones multiescala
(ver \cite{Cohen}). Llamaremos {\em raíz} a la base, y \emph{nodo}
a cada elemento del árbol. Un nodo padre posee cuatro hijos, y un nodo
sin hijos es llamado \emph{hoja}. Cada nodo posee $s'=2$ vecinos en cada
 dirección espacial, llamados \emph{primos cercanos}, necesarios para
 determinar los flujos en cada hoja; si tales primos cercanos no existieran,
 deben ser creados artificialmente como \emph{hojas virtuales}. Las hojas
 del árbol son los elementos que conforman la malla adaptativa. Denotamos por
 $\mathcal{L}(\Lambda)$ a la restricción del conjunto de nodos $\Lambda$ al
conjunto de hojas.

En cada paso temporal, la reconstrucción de multiresolución es aplicada
a la parte espacial de la solución $\u=(v,u_\mathrm{e},w)$. Luego de cada
paso temporal, es necesario actualizar la malla adaptativa, y esto se realiza
mediante la aplicación de la estrategia de corte. Una vez que se aplica tal
estrategia, es necesario agregar una \emph{zona de seguridad} a la estructura
de árbol con el fin de asegurar que la malla a ser utilizada en el paso temporal
siguiente, será adecuada para representar la solución correspondiente.
Esta zona de seguridad será implementada agregando un nivel más fino en todas las
posiciones posible (sin destruir la estructura de árbol).
Para cuantificar la mejora obtenida en compresión de datos y tiempo de máquina,
usaremos la
\emph{tasa de compresión de datos} y la \emph{aceleración de tiempo de máquina}
 (ver \cite{brss2})
\begin{equation*}
\eta:=\frac{\mathcal{N}}{2^{-(L+1)}\mathcal{N}+\#\mathcal{L}(\Lambda)},\quad
\mathcal{V}:=\frac{\mathrm{CPU\, time}_{\mathrm{FV}}}{\mathrm{CPU\, time}_{\mathrm{MR}}}.
\end{equation*}
Aquí  $\mathcal{N}$ es el número de elementos en la malla más fina en el nivel $L$, y
$\#\mathcal{L}(\Lambda)$ es la cardinalidad del conjunto de hojas.

\subsection{Análisis de error para el método de multiresolución}
\label{sec:err-analysis}
Usando las propiedades básicas del esquema de volúmenes finitos
de referencia, derivamos la elección óptima para la tolerancia de
referencia (\ref{eq4.4}). En primer lugar, el error global entre
la solución de referencia y la solución mediante multiresolución
es descompuesto en dos errores
$$\bigl\|\u^L_{\mathrm{ex}} - \u^L_{\mathrm{MR}} \bigr\| \leq
\bigl\|\u^L_{\mathrm{ex}} - \u^L_{\mathrm{FV}} \bigr\| +
\bigl\|\u^L_{\mathrm{FV}} - \u^L_{\mathrm{MR}} \bigr\|.$$
El primer error del lado derecho es denominado \emph{error de
discretización} y el segundo error es denominado \emph{error de
perturbación}. Utilizando estimaciones estándar para ambos errores
y la condición CFL (\ref{cfl}), obtenemos (ver detalles en \cite{brss,Cohen})
que si la \emph{tolerancia de referencia} es dada por
\begin{equation*}\epsilon_{\mathrm{R}} = C \frac{2^{(2-\alpha)L-2}}{\displaystyle{
\max_{K\in\TT}\Bigl(|\IonK|+2|\IapK|\Bigr)+D\,
\max_{K\in\TT}\Bigl(|M_{i,K}|+|M_{e,K}|\Bigr)}},
\end{equation*}
entonces el error de discretización y el error de perturbación poseen el mismo
orden de magnitud.

\begin{figure}[ht]
\begin{center}
\begin{tabular}{ccc}
\includegraphics[width=0.3\textwidth]{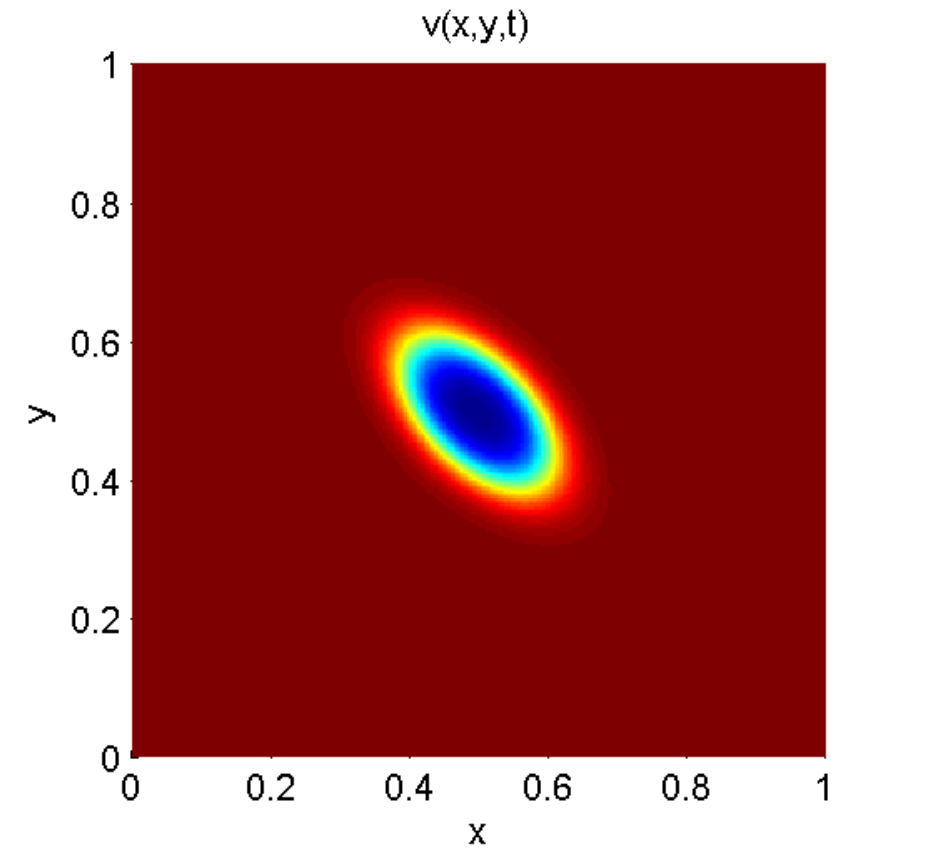}&
\includegraphics[width=0.3\textwidth]{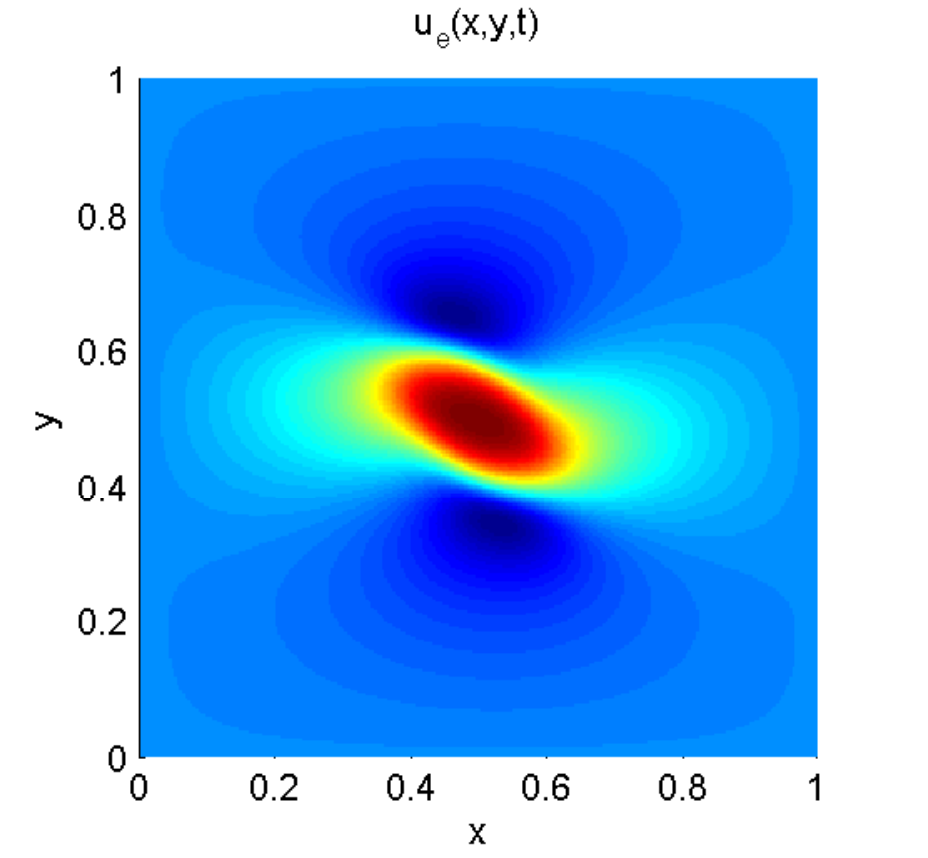}&
\includegraphics[width=0.3\textwidth]{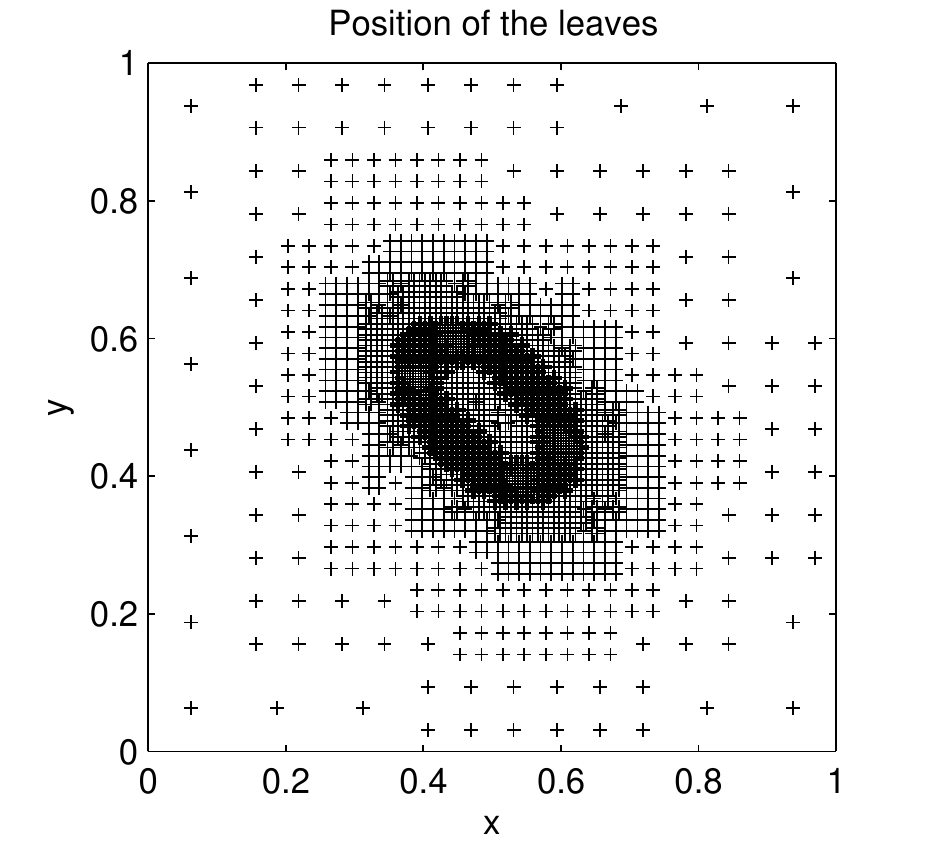}\\
\includegraphics[width=0.3\textwidth]{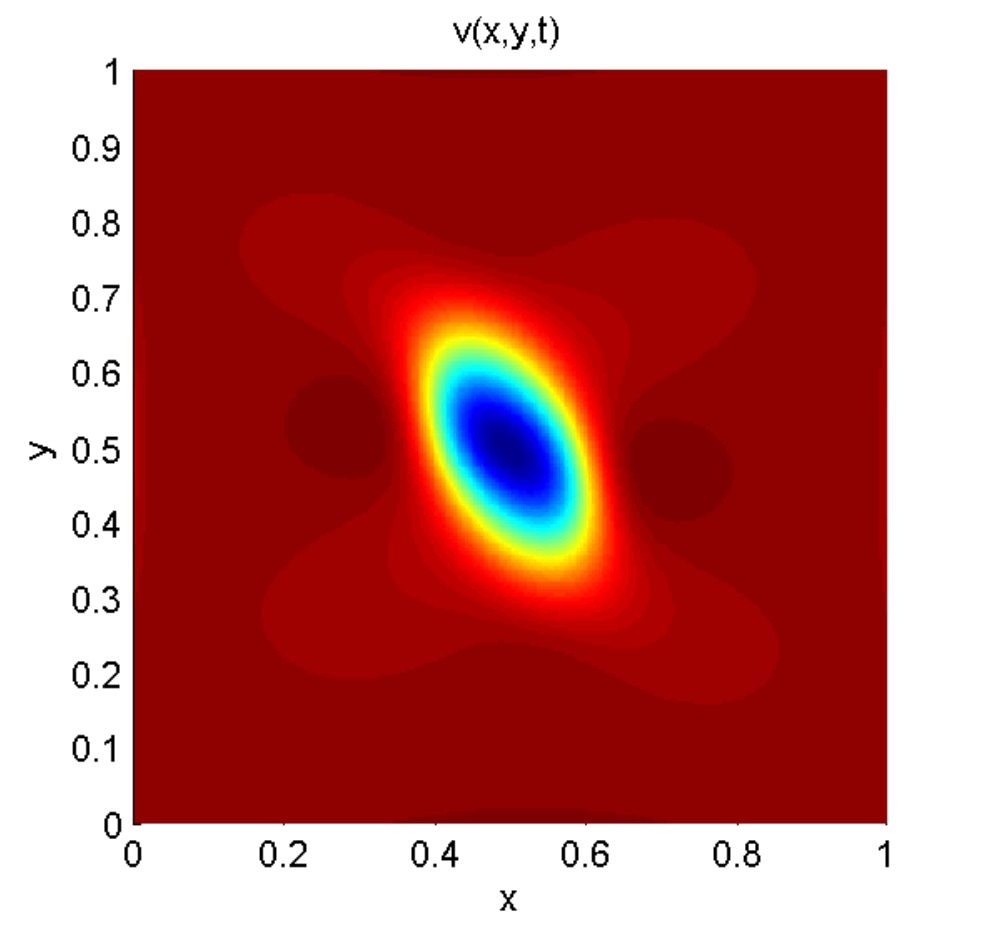}&
\includegraphics[width=0.3\textwidth]{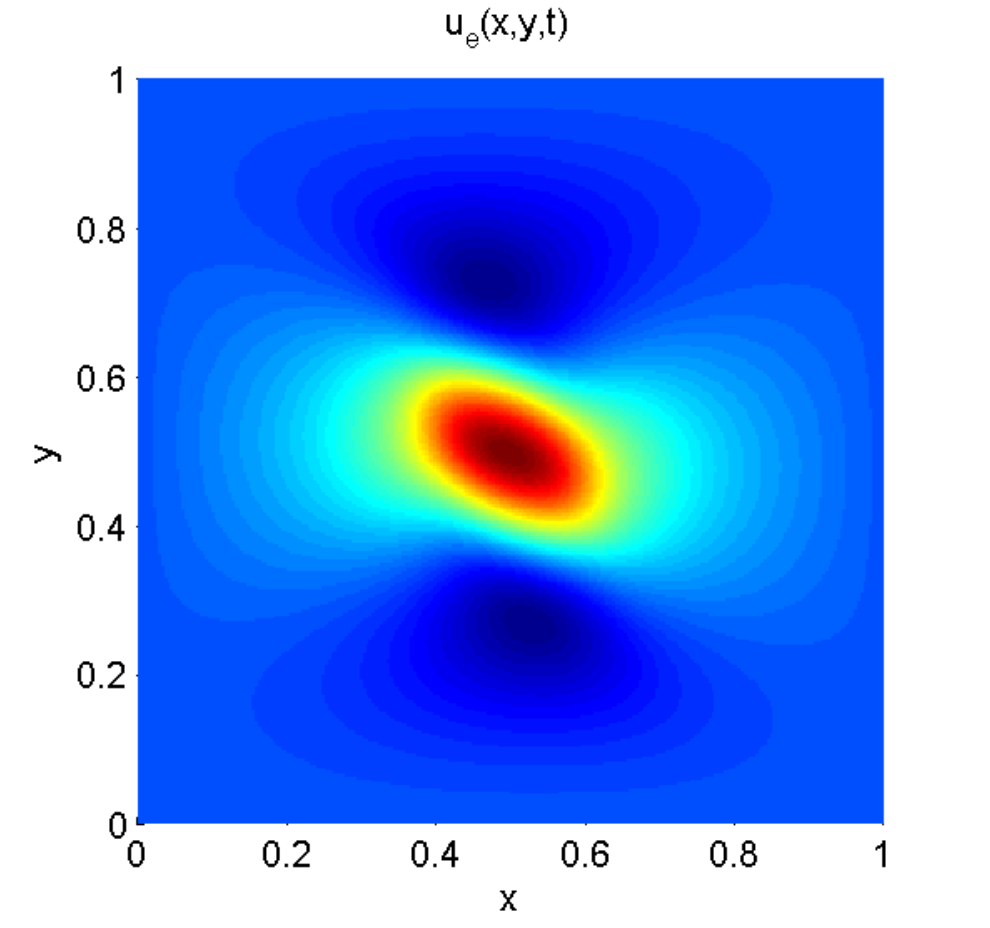}&
\includegraphics[width=0.3\textwidth]{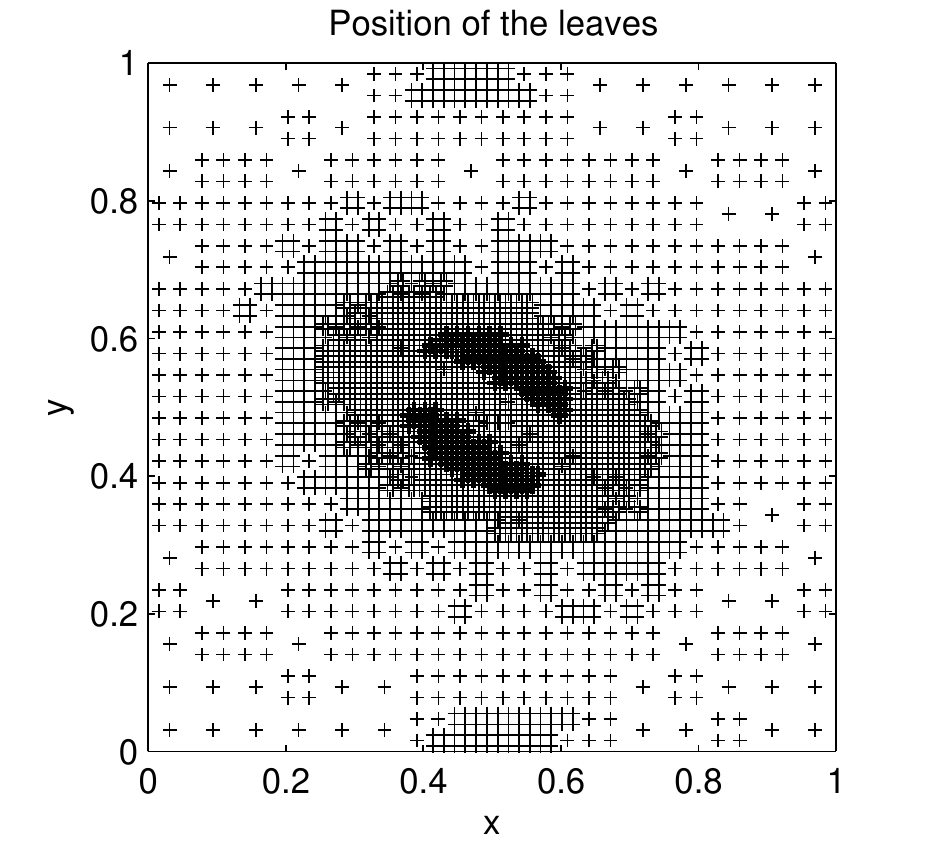}
\end{tabular}
\caption{\it Potencial transmembrana $v$,  potencial
extracelular $u_\mathrm{e}$,  y hojas de la correspondiente estructura de árbol
 en los instantes $t=0.01$ y $t=1.1$.} \label{fig:snapshots1}
\end{center}
\end{figure}


\subsection{Aceleración de la evolución temporal}
La idea básica del método a ser presentado, es utilizar una condición CFL
local, imponiendo el mismo número CFL para todas las escalas. La estrategia
consiste en evolucionar todas la hojas situadas en el nivel $l$ usando el
paso temporal local
\begin{equation*}
\Delta t_l=2^{L-l}\Delta t,\quad l=L-1,\ldots,0,
\end{equation*}
donde $\Delta t=\Delta t_L$ corresponde al paso temporal sobre el nivel más
fino $L$. Tal estrategia permite incrementar el paso temporal para la mayor
parte de la malla adaptativa, sin necesidad de violar la condición CFL. Ahora,
es necesario sincronizar el paso temporal para las porciones de la solución
que se encuentran en distintos niveles de resolución. Pero esta sincronización
es alcanzada de manera automática al cabo de $2^l$ pasos usando  $\Delta t_l$.

La proyección y predicción de multiresolución son efectuadas sólo en los niveles
ocupados por hojas del árbol correspondiente, y sólo cada dos pasos temporales.
Para el resto de los pasos intermedios, se utiliza la misma estructura de árbol.
Del mismo modo, los flujos son calculados sólo en niveles que contienen hojas,
y éstos son calculados del siguiente modo: Si un borde dado es compartido por
dos hojas en el mismo nivel $l$, entonces el cálculo del flujo se hace de manera
estándar, utilizando los primos cercanos u hojas virtuales si fuere necesario.
Si el borde es compartido por una hoja en el nivel $l$ y una hoja en un nivel
más fino $l+1$ (\emph{borde interfaz}), calculamos los flujos en el nivel $l+1$
sobre el mismo borde, y el flujo correspondiente en el nivel $l$ será igual a la
suma de los flujos (en la dirección opuesta) sobre los hijos correspondientes en
el nivel $l+1$. Con el fin de tener siempre a disposición los flujos calculados,
la estrategia de paso temporal local debe realizarse recursivamente desde el nivel
más fino hasta el nivel más grueso.
{\small
\begin{table}[htb]
\begin{center}
\begin{tabular}{lccccccc}
\hline
Tiempo $[\mathrm{ms}]$& $\mathcal{V}$  & $\eta$&  Potencial & error $L^1$
& error $L^2$ &error $L^\infty$ $\vphantom{\int^X}$   \\
\hline
$t=$0.01  $\vphantom{\int^X}$
& 13.74 & 19.39& $v$  &$3.68\times10^{-4}$
&$8.79\times10^{-5}$&$6.51\times10^{-4}$ \\
&      &      & $u_\mathrm{e}$ &$2.01\times10^{-4}$
&$6.54\times10^{-5}$&$5.22\times10^{-4}$  \\
$t=$1.1
& 21.40 & 17.63 & $v$& $4.06\times10^{-4}$
&$9.26\times10^{-5}$&$6.83\times10^{-4}$ \\
&      &       & $u_\mathrm{e}$&$2.79\times10^{-4}$
&$8.72\times10^{-5}$&$5.49\times10^{-4}$ \\
$t=$2.2
& 25.23 & 17.74& $v$ & $4.37\times10^{-4}$
&$1.25\times10^{-4}$&$6.88\times10^{-4}$ \\
&      &     & $u_\mathrm{e}$ & $3.48\times10^{-4}$
&$9.44\times10^{-5}$&$6.11\times10^{-4}$ \\
$t=$3.3
& 26.09 & 16.35 & $v$& $5.29\times10^{-4}$
&$1.94\times10^{-4}$&$7.20\times10^{-4}$ \\
&       &     & $u_\mathrm{e}$& $4.15\times10^{-4}$
&$1.06\times10^{-4}$&$6.32\times10^{-4}$ \\
\hline
\end{tabular}
\end{center}
\caption{\it Tiempo de simulación, aceleración de tiempo de máquina~$\mathcal{V}$, tasa
de compresión~$\eta$ y errores normalizados.}
\label{table:ex2}
\end{table}}

\section{Ejemplo numérico}
\begin{figure}[ht]
\begin{center}
\begin{tabular}{ccc}
\includegraphics[width=0.3\textwidth]{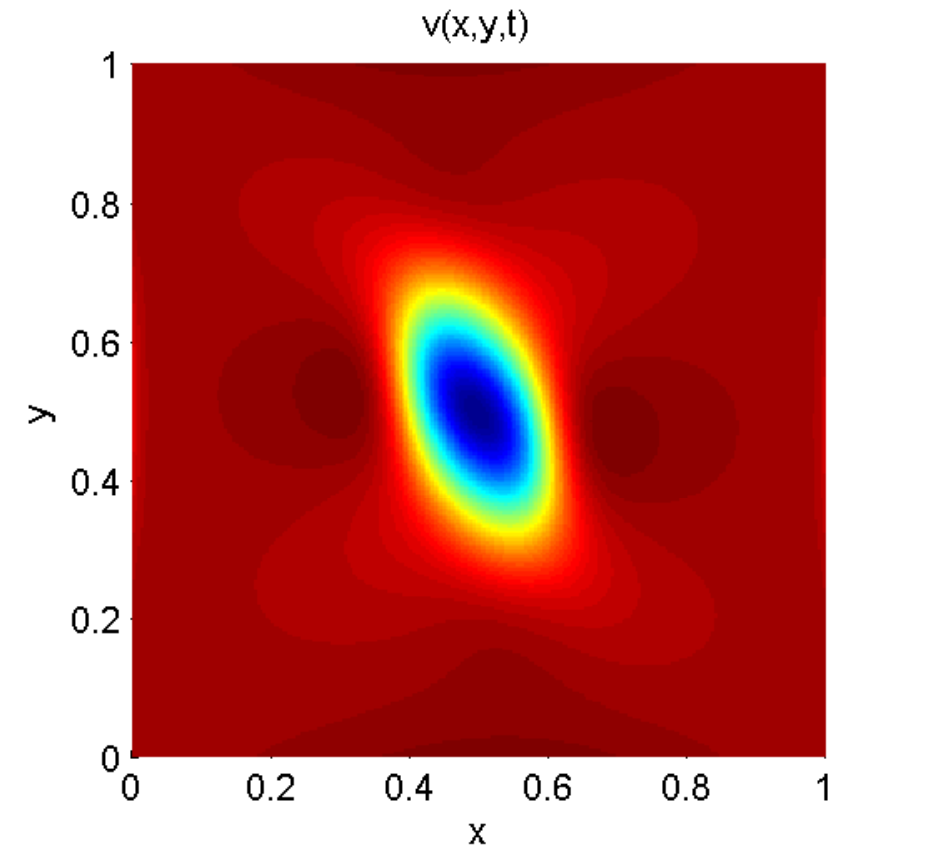}&
\includegraphics[width=0.3\textwidth]{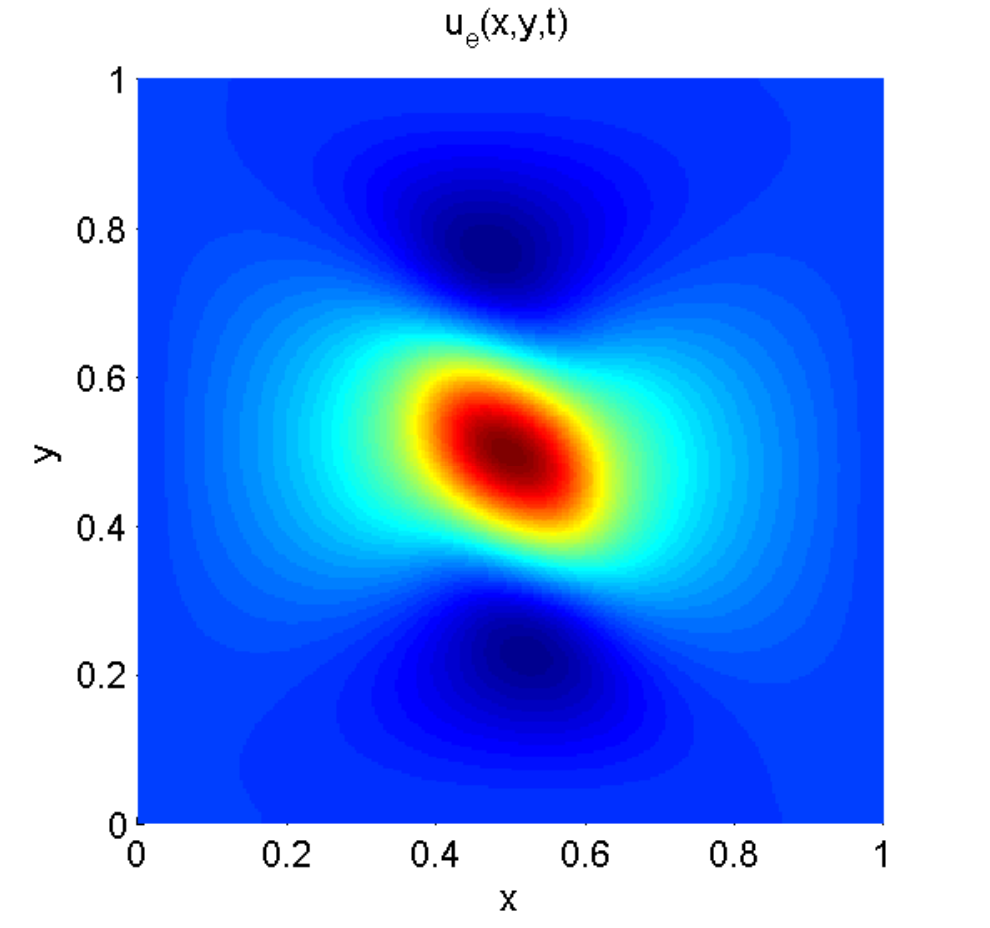}&
\includegraphics[width=0.3\textwidth]{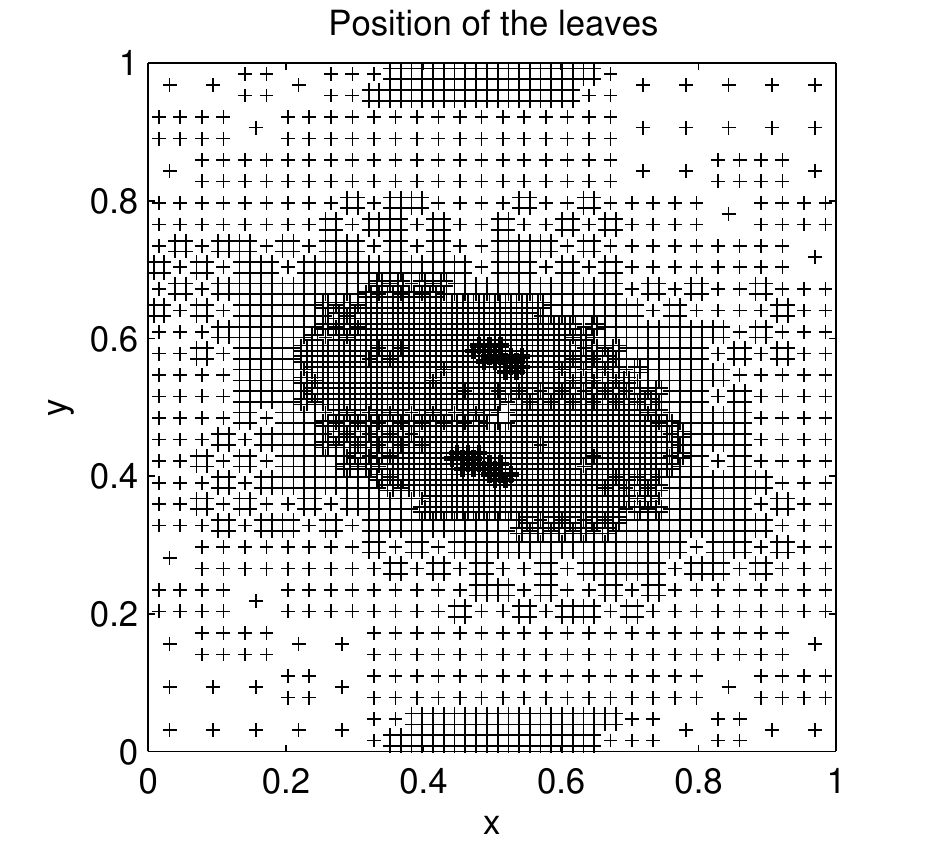} \\
\includegraphics[width=0.3\textwidth]{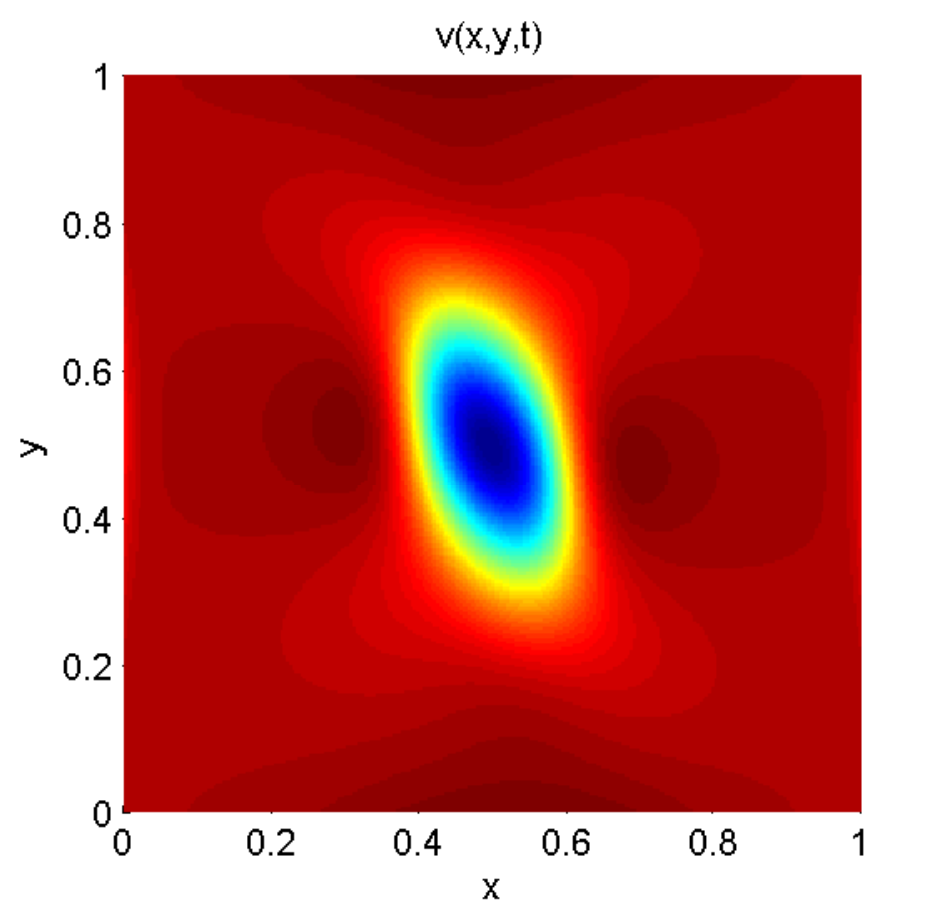}&
\includegraphics[width=0.3\textwidth]{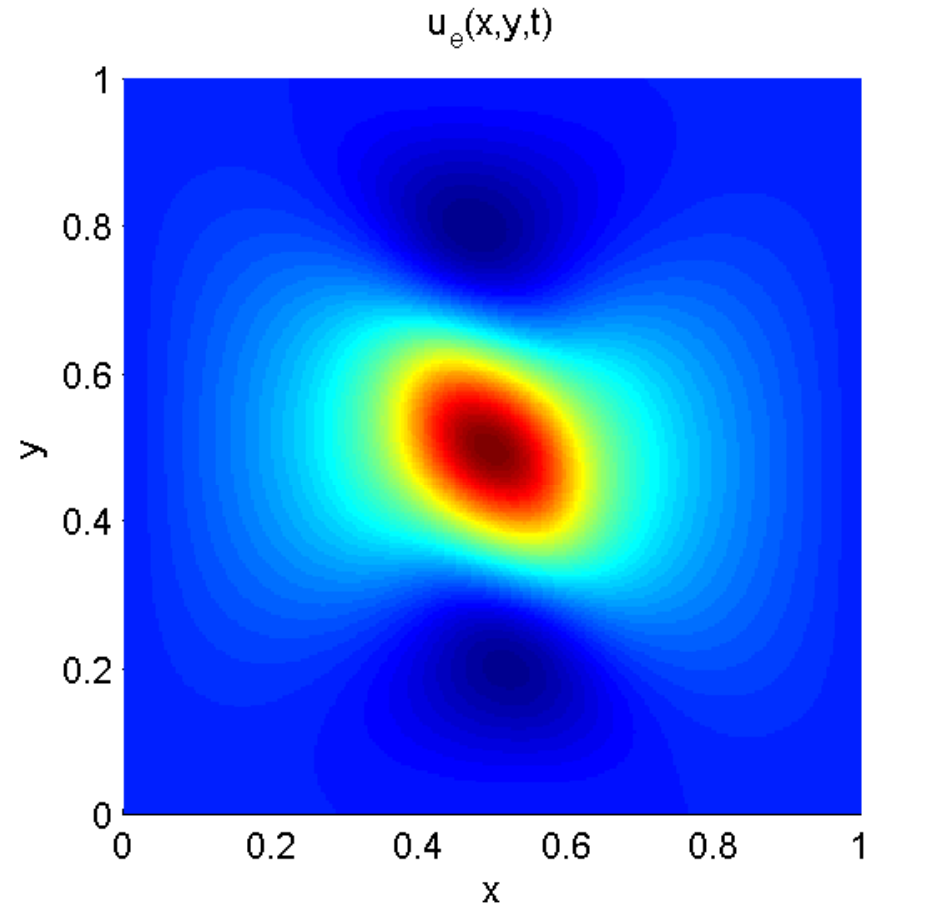}&
\includegraphics[width=0.3\textwidth]{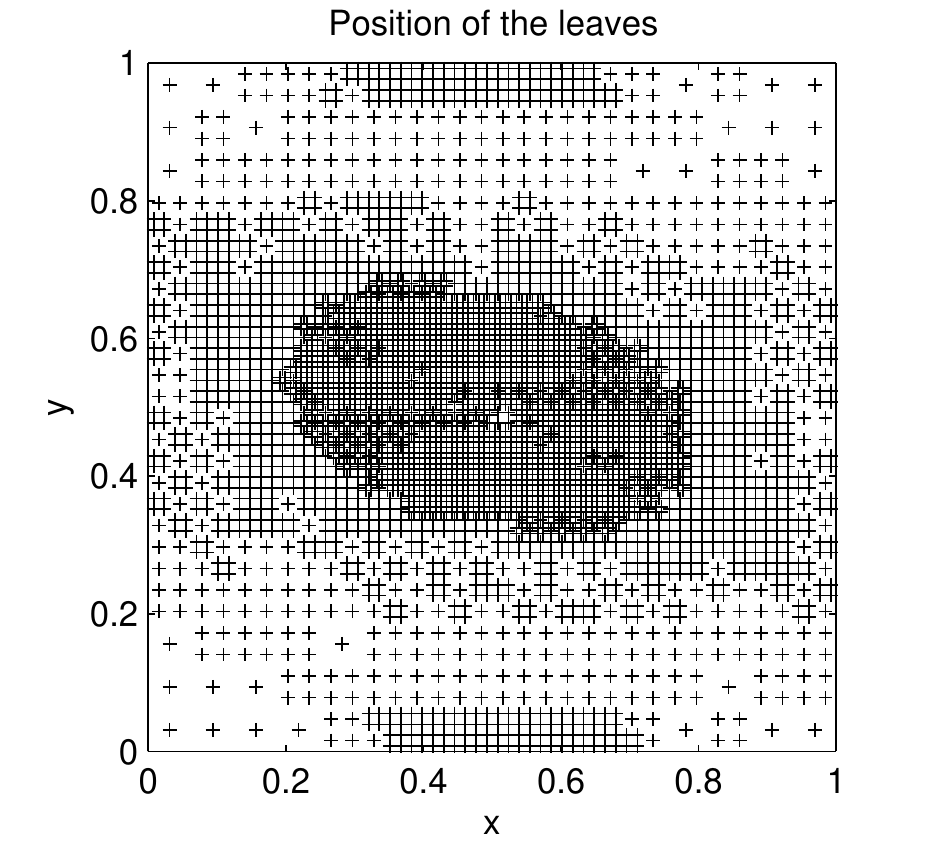}
\end{tabular}
\caption{\it Potencial transmembrana $v$,  potencial
extracelular $u_\mathrm{e}$,  y hojas de la correspondiente estructura de árbol
 en los instantes $t=2.2$ y $t=3.3$.} \label{fig:snapshots2}
\end{center}
\end{figure}

En las simulaciones se utiliza un dominio computacional simple
$\Om=[0,1\,\mathrm{cm}]^2$ y los siguientes parámetros (siguiendo \cite{ying}):
capacitancia de la membrana $c_m=1.0\,\mathrm{mF/cm}^2$, conductividades
$\sigma_\mathrm{i}^l=6\,\Om^{-1}\mathrm{cm}^{-1}$,
$\sigma_\mathrm{i}^t=0.6\,\Om^{-1}\mathrm{cm}^{-1}$,
 $\sigma_\mathrm{e}^l=24\,\Om^{-1}\mathrm{cm}^{-1}$
y $\sigma_\mathrm{e}^t=12\,\Om^{-1}\mathrm{cm}^{-1}$, razón superficie volumen
 $\beta=4036.5\,\mathrm{cm}^{-1}$, resistividad de superficie
$R_m=2\times10^{4}\,\Om\,\mathrm{cm}^2$, $v_p=100\,\mathrm{mV}$, $\eta_1=0.005$,
$\eta_2=0.1$, $\eta_3=1.5$, $\eta_4=7.5$, $\eta_5=0.1$. Las fibras forman un
ángulo de $-\pi/4$ con el eje $x$ y como dato inicial, aplicamos un estímulo
en el medio extracelular en el centro del dominio (ver
Figura~\ref{fig:snapshots1}). Se elige la siguiente configuración para el método de
multiresolución: Wavelets con $r=3$ momentos nulos, nivel maximal de resolución
$L=9$ y por lo tanto una malla fina de $\mathcal{N}=65536$ elementos, una
tolerancia de referencia dada por $\eps_R=5.0\times 10^{-4}$. Mostramos en las
Figuras~\ref{fig:snapshots1},\ref{fig:snapshots2} una secuencia de instantáneas de
la solución después de haber aplicado un estímulo en el centro del dominio. Se
muestra tanto la solución, como la correspondiente malla adaptativa generada por la
multiresolución. Los errores han sido calculados utilizando como referencia, una
solución aproximada de volúmenes finitos sobre una malla con
$\mathcal{N}=1024^2=1048576$ volúmenes de control. En la tabla \ref{table:ex2}
puede notarse que la solución numérica obtenida aplicando multiresolución es
suficientemente precisa (errores del orden de $10^{-4}$) y las tasas de compresión
son considerablemente altas.

Para la integración temporal usando LTS, elegimos $\mathrm{CFL}_0=0.5$ para
el nivel más grueso de resolución, y $\mathrm{CFL}_l=2^l\mathrm{CFL}_0$ para los
niveles más finos. Utilizando LTS, se obtiene un aumento sustancial
en tasa de aceleración, y sin embargo los errores se mantienen con el mismo orden
de precisión.

\section{Agradecimiento}
MB agradece el apoyo del proyecto Fondecyt 1070682, RB agradece
el apoyo del proyecto Fondecyt 1050728 y el programa Fondap en Matemática Aplicada,
 proyecto 15000001; y RR agradece el apoyo de Beca Conicyt.


\end{document}